\documentstyle[12pt]{article} 

\title{Weierstrass representations for harmonic morphisms on 
Euclidean spaces and spheres}

\author{P. Baird and J. C. Wood
\thanks{Partially supported by EC grant CHRX-CT92-0050} \\
{\small D\'epartement de Math\'ematiques, Universit\'e de Bretagne
Occidentale} \\
{\small 6 Avenue Le Gorgeu, B.P. 452, 29275 Brest C\'edex, France,
and} \\
{\small Department of Pure Mathematics, University of Leeds}\\
{\small Leeds LS2 9JT, G.B.} }
\date{}

\newtheorem{theo+}           {Theorem}      [section]
\newtheorem{prop+}  [theo+]  {Proposition}
\newtheorem{coro+}  [theo+]  {Corollary}
\newtheorem{lemm+}  [theo+]  {Lemma}
\newtheorem{exam+}  [theo+]  {Example}
\newtheorem{rema+}  [theo+]  {Remark}
\newtheorem{defi+}  [theo+]  {Definition}

\newtheorem{exam+s}  [theo+]  {Examples}
\newtheorem{rema+s}  [theo+]  {Remarks}
\newtheorem{hyp+}  [theo+]  {Hypotheses}
\newtheorem{cla+}  [theo+]  {Claim}

\newenvironment{theorem}{\begin{theo+}}{\end{theo+}}
\newenvironment{proposition}{\begin{prop+}}{\end{prop+}}

\newenvironment{lemma}{\begin{lemm+}}{\end{lemm+}}
\newenvironment{example}{\begin{exam+}\rm}{\end{exam+}}

\newenvironment{definition}{\begin{defi+}\rm}{\end{defi+}}

\newcommand{\Bbb}{\bf} 
 
\newcommand{\RR}{{\Bbb R}} 
\newcommand{\CC}{{\Bbb C}}  
\newcommand{\NN}{{\Bbb N}}  
\newcommand{\Rn}{{\Bbb R}^n} 
\newcommand{\Rmm}{{\Bbb R}^{2m}} 
\newcommand{\Cm}{{\Bbb C}^m} 

\newcommand{\ii}{{\rm i}}  

\newcommand{\proof}{\noindent\ {\bf Proof} \hskip 0.4em}
\newcommand{\eproof}{\bigskip}

\begin{document} 
\maketitle

\begin{abstract}
We construct large families of harmonic morphisms which are
holomorphic with respect to Hermitian structures by finding
heierarchies of Weierstrass-type representations. This enables us
to find new examples of complex-valued harmonic morphisms from Euclidean
spaces and spheres.  
\end{abstract}

\section{Introduction}\label{sec:Introduction}
Let $\phi :U\rightarrow N$ be a smooth mapping between Riemannian
manifolds. Then $\phi$ is called a {\em harmonic morphism\/} if,
for every real valued function harmonic on an open set $W\subset
N$ with $\phi ^{-1}(W)$ non-empty, the pull-back $f\circ \phi$ is
harmonic on $\phi ^{-1}(W)$ in $M$. Equivalently, $\phi$ is a
harmonic morphism if and only if $\phi$ is both horizontally weakly
conformal and harmonic \cite{Fu-1,Is}.  In the case when $U$ is
a domain in $\Rn$ with its standard Euclidean structure and
$N$ is the complex plane $\CC$, the equations for horizontal weak
conformality and harmonicity are, respectively,
\begin{equation}
\sum_{i=1}^n\left( \frac{\partial\phi}{\partial x^i}\right) ^2
=0\;,
\label{eq:one}
\end{equation}
\begin{equation}
\triangle\phi\equiv\sum_{i=1}^n\frac{\partial ^2\phi}
{(\partial x^i)^2}=0\;.\label{eq:two}
\end{equation}
In fact, provided (\ref{eq:one}) is satisfied, (\ref{eq:two})
is equivalent to the minimality of the fibres at regular points
\cite{Ba-Ee}.

An important example of a harmonic morphism is the following.  Let
$U$ be a domain of ${\RR}^{2m}\cong \CC^m$ and suppose that
$\phi$ is holomorphic with respect to a K\"ahler structure on
${\RR}^{2m}$, then $\phi$ is a harmonic morphism \cite{Fu-1}. 
In a more general setting, but restricting the dimension of the
domain to be 4, the second author showed that, if $U$ is a domain
of a 4-dimensional Einstein manifold, then any submersive
harmonic morphism $\phi$
is holomorphic with respect to an integrable Hermitian
structure $J$ on $U$.  Furthermore, the fibres of
$\phi$ are superminimal with repect to $J$ \cite{Wo}. 
Thus there is a strong relationship between holomorphic maps and
harmonic morphisms.

In \cite{Ba-Wo-5}, we studied those harmonic morphisms defined on
domains $U$ of ${\Rmm}, (m\in {\NN})$, which are holomorphic with respect
to an integrable Hermitian structure on $U$, finding full global
examples in the case $m\geq 3$, which are neither holomorphic with respect
to a K\"ahler structure nor have superminimal fibres, as well as
such examples which factor to a domain of ${\RR}^{2m-1}$. In
this paper we concentrate on the case of superminimal fibres and
study reduction to odd-dimensional Euclidean spaces and to spheres,
constructing large classes of examples in terms of holomorphic
data.

More generally, we refer to a parametrization of solutions to
Equations (\ref{eq:one}) and (\ref{eq:two}) in terms of
holomorphic data as a ``Weierstrass representation", after the
local representation of Weierstrass for minimal surfaces in ${\RR}^3$.
We have observed an interesting duality between the
theory of minimal surfaces and harmonic morphisms
\cite{Ba-Wo-1}.

In \cite{Ba-Wo-1}, Weierstrass representations for harmonic
morphisms defined on Euclidean spaces and spheres were obtained
in the case when the fibres are totally geodesic.  More
precisely, if $(\xi _1, \ldots ,\xi _n)$ is an $n$-tuple of
meromorphic functions of $z$ such that $\sum \xi _i^2=0$, then:

(i) the {\em inhomogeneous\/} equation 
\begin{equation}
\xi _1x^1+\ldots +\xi _nx^n=1
\label{eq:tginhomog}
\end{equation}
locally determines all (submersive complex-valued) harmonic morphisms
$z=z(x)$ on domains of $\Rn$ with totally geodesic fibres;

(ii)  the {\em homogeneous\/} equation 
\begin{equation}
\xi _1x^1+\ldots +\xi _nx^n=0
\label{eq:tghomog}
\end{equation}
locally determines all (submersive complex-valued) harmonic morphisms
$z=z(x)$ on domains of $S^{n-1}$ with totally geodesic fibres.  A similar
representation
for harmonic morphisms from real hyperbolic spaces $H^{n-1}$ with totally
geodesic fibres has been obtained by S. Gudmundsson \cite{Gu-3}.

Here we develop a far richer description of harmonic morphisms in
terms of holomorphic data, obtaining a {\em heierarchy\/} of
Weierstrass representations on domains of ${\RR}^n$ such that:

(i) on ${\Rmm}$, the Weierstrass representation locally gives all
(submersive complex-valued) harmonic
morphisms which are holomorphic with respect to a Hermitian structure with
superminimal fibres;

(ii)  reduction ${\RR}^n\rightarrow {\RR}^{n-1}$ commutes with
the Weierstrass representation;

(iii)  on ${\RR}^3$ and ${\RR}^4$, the Weierstrass representation
describes {\em all\/} (submersive complex-valued) harmonic morphisms.

Finally we show how to construct harmonic morphisms from spheres
using our Weierstrass representations by choosing the appropriate
holomorphic function to be homogeneous. We refer the reader to
the work of Gudmundsson for other constructions of harmonic
morphisms on complex projective spaces and K\"ahler manifolds
\cite{Gu-1,Gu-2}.

\bigskip

The authors would like to thank the referee for helpful comments on this
work.

\section{Holomorphic harmonic morphisms and reduction}
\label{sec:red}
By a Hermitian structure on an open subset $U$ of $\Rmm$ we mean
the smooth choice of an almost Hermitian structure on
$U$ which is integrable (cf. \cite{Ba-Wo-5}).  By a {\em
holomorphic harmonic morphism\/} on $U$ we mean a complex-valued
harmonic morphism which is holomorphic with respect to some
Hermitian structure on $U$. We briefly summarize the
characterization of holomorphic harmonic morphisms given in
\cite{Ba-Wo-5} and describe how to reduce to other manifolds.

Let $(x^1,\ldots ,x^{2m})$ be standard coordinates for ${\Rmm}$
and introduce complex coordinates $q^1=x^1+ix^2, q^2=x^3+ix^4$
etc. If $J$ is a positive integrable Hermitian structure defined
on an open subset $U\subset {\Rmm}$, then locally $J$ is
characterised by $m(m-1)/2$ functions $\mu _1, \ldots ,\mu
_{m(m-1)/2}$ holomorphic in $m$ complex parameters $(z^1,
\ldots , z^m)$ as follows:

Given $\mu = (\mu _1, \ldots ,\mu _{m(m-1)/2})$, let $M=\left(
M^i_{\bar j}(\mu)\right) $ be the skew symmetric matrix
\begin{equation}
\big(M^i_{{}{\bar{j}}}(\mu)\big) =
\left(\begin{array}{rrrrr} 0          & {}\mu_1    & {}\mu_2    &
\ldots & \mu_{m-1} \\ -\mu_1    & 0           & {}\mu_m  & \ldots
& \mu_{2m-3} \\ -\mu_2     & -\mu_m      & 0         & \ldots &
\mu_{3m-6} \\ \vdots     & \vdots      & \vdots      & \ddots &
\vdots   \\ -\mu_{m-1} & -\mu_{2m-3} & -\mu_{3m-6}  & \ldots & 0
\end{array}\right)\;. \label{M} 
\end{equation}
As in \cite{Ba-Wo-5} this matrix determines a positive almost
Hermitian structure at any point of $\Rmm$, namely that with
(1,0)-cotangent space given by $\mbox{span}
\{e^i=dq^i-M^i_{\bar{j}}dq^{\bar{j}}: i=1, \ldots , m\}$.~(\footnote{We
use the double summation convention throughout.})
Now
suppose that the $\mu _i = \mu _i(z)$ are holomorphic functions on a domain
$V$  of $\Cm$ and let $h^1(z), \ldots , h^m(z)$ be
further holomorphic functions on $V$.  Consider the following system of
equations:
\begin{equation}
F^i(q,z)\equiv q^i-M^i_{\bar{j}}(z)q^{\bar{j}}-h^i(z)=0,
 \quad (i=1,\ldots,m) \;.
\label{eq:three}
\end{equation}
where we write $M(z)$ for $M(\mu(z))$.
The system (\ref{eq:three}) has the form 
\begin{equation}
F(q,z)=0
\label{eq:four}
\end{equation}
where $F:\Rmm\times V  \rightarrow {\Cm}$.  Provided the
determinant of the  Jacobian matrix 
$$
K=\left( \partial _jF^i\right) \mbox{ where }
\partial _j=\frac{\partial}{\partial z^j}
$$
is non-zero, we can locally solve (\ref{eq:four}) for $z=z(q)$. 
Then, on suitable open sets, $z(q)=(z^1(q), \ldots ,z^m(q))$ form complex
coordinates for the complex manifold $(U, J)$ where $J=J(M)$ is
given at $q \in U$ by $J_q=J(M(z(q)))$.  Note that all complex
coordinates with respect to any Hermitian structure are given locally
this way. Indeed, if $q\mapsto M^i_{\bar{j}}(q)$ defines a
Hermitian structure, then
\begin{equation}
w^i=q^i-M^i_{\bar{j}}(q)q^{\bar{j}} \quad (i=1,\ldots,m)
\label{eq:wi}
\end{equation}
give complex coordinates 
(see, for example \cite{Ba-Wo-5}).  Any other complex coordinates 
$z=(z^1, \ldots , z^m)$ are related by equations $w^i=h^i(z)$ for some 
holomorphic functions $h^i$, whence, writing
$M^i_{\bar{j}}(q)=M^i_{\bar{j}}(z(q))$, we see that $z=z(q)$ satisfies
(\ref{eq:three}). 

More precisely, only those Hermitian structures with values at each point
in a
large cell of the space $SO(2m)/U(m)$ of all positive almost Hermitian
structures are parametrized as above.  However, by acting on
${\Rmm}$ by an isometry, we may always assume, at least locally,
that $J$ satisfies this condition.  In the general case, we must
allow the $\mu _i$'s to become infinite.

Let $J$ be a Hermitian structure on $U\subset {\Rmm}$
characterized by Equation (\ref{eq:three}) above.  Then any
holomorphic map $\phi :(U, J) \rightarrow \CC$ is a
holomorphic function of the complex coordinates $(z^1, \ldots
,z^m)$.  Now, in a neighbourhood of a regular point of $\phi$, it is no
loss of
generality to suppose that $\phi =z^1$.  Then the Laplacian of the
function $z^1=z^1(q)$ is calculated to be \cite{Ba-Wo-5}:
\begin{equation}
\Delta z^1 = \frac{4}{\det K} \sum_{1 \leq i<j \leq m}
(-1)^{i+j}
\left| \begin{array}{cccc}
\partial _2M^i_{\bar{j}} & \partial _3M^i_{\bar{j}} & \ldots
 & \partial _mM^i_{\bar{j}} \\
\partial _2F^{k_1} & \partial _3F^{k_1} &  \ldots &
\partial_mF^{k_1} \\
\vdots & \vdots & \ddots & \vdots \\
\partial _2F^{k_{m-2}} & \partial _3F^{k_{m-2}} &  \ldots &
\partial _mF^{k_{m-2}} 
\end{array}
\right| \label{eq:five}
\end{equation}
where $(k_1, \ldots , k_{m-2})=(1, \ldots , \widehat{i}, \ldots ,
 \widehat{j}, \ldots , m)$.  
Since $z^1$ is holomorphic it is horizontally weakly conformal, and so, if
$\triangle z^1 =0$ then $z^1:U\rightarrow \CC$ is a harmonic
morphism.

A particular class of solutions is given by the case when each
$\mu _i(z)$ depends on $z^1$ only, for then all the partial
derivatives $\partial _2M^i_{\bar{j}}, \ldots , \partial
_mM^i_{\bar{j}}$ vanish. In this case the Hermitian structure $J$
is constant along the fibres of $z^1$, equivalently the
fibres of $z^1$ are {\em superminimal}. Conversely, by the
above remarks, any (submersive complex-valued) harmonic morphism defined
on an open subset of
${\Rmm}$, holomorphic with respect to a Hermitian structure and
with superminimal fibres, may be described this way locally.

Our method of reducing holomorphic harmonic morphisms defined on
domains in ${\Rmm}$ to other manifolds follows from the
fundamental property that the composition of two harmonic
morphisms $\phi \circ \pi$
is a harmonic morphism \cite{Fu-1} and a converse to
this \cite[Proposition 1.1]{Gu-1} (note that the author is assuming
surjectivity of $\pi$
without mentioning it).  In particular, we shall be concerned with the
following projections which are all harmonic morphisms:

(i)  {\em Orthogonal projection\/}  $\pi ^n_1:\Rn \rightarrow
{\RR}^{n-1}$ given by the formula $\pi
^n_1(x^1, \ldots , x^n)=(x^2, \ldots , x^n)$.

(ii)  {\em Radial projection\/}  $\pi ^n_r: \Rn\backslash
\{ 0\}\rightarrow S^{n-1}$ given by the formula $\pi ^n_r(x^1,
\ldots , x^n)=(x^1, \ldots , x^n)/|x|$.

(iii) {\em Radial projection followed by the Hopf map\/}  $\pi
^{2m}_H: {\Rmm}\backslash \{ 0\} = \Cm\backslash \{ 0\}
\rightarrow {\CC}P^{m-1}$ given
by the formula $\pi ^{2m}_H(q^1, \ldots , q^m) = [q^1, \ldots ,
q^m]$.

So, for example, a harmonic morphism on a domain in $S^{n-1}$ is
equivalent to a harmonic morphism on a domain in $\Rn
\backslash \{ 0\}$ invariant under radial projection. 
Precisely, if $\phi :W\subset S^{n-1} \rightarrow {\CC}$ is a
harmonic morphism, then the composition ${ \Phi =\phi\circ\pi
^n_r|_{(\pi ^n_r)^{-1}(W)}}$ is a harmonic morphism on $(\pi
^n_r)^{-1}(W)$ such that $\partial\Phi/\partial r=0$, where
$r=|x|$ denotes the radial coordinate.  Conversely, if $\Phi :
U\subset \Rn \rightarrow {\CC}$ is a harmonic morphism on
an open set $U$ with $U\cap S^{n-1}$ non-empty and such that
$\partial\Phi/\partial r=0$, then it follows that the restriction
$\phi = \Phi |_{U\cap S^{n-1}}:U\cap S^{n-1} \rightarrow {\CC}$
is a harmonic morphism on the domain $U\cap S^{n-1} \subset
S^{n-1}$.  We formulate the invariance in a more general
setting:

Let $J$ be a Hermitian structure defined on a domain $U\subset
{\Rmm}$ with associated complex coordinates $(z^1, \ldots , z^m)$
given locally by Equation (\ref{eq:three}).  Let
$v=a^j(x) (\partial/ \partial x^j)$ be a vector field on
$U$.  If $\phi : (U, J)\rightarrow {\CC}$ is a holomorphic
function, then $\phi$ is {\em invariant under\/} $v$ if and only
if the directional derivative $d\phi (v)$ vanishes.  As remarked
above, in a neighbourhood of a regular point of $\phi$
we may assume without loss of generality that $\phi$ is the
holomorphic function $z^1$.

\begin{proposition}\label{pr:red} Let
$\phi =z^1:(U,J)\rightarrow {\CC}$ be a holomorphic function,
where $z = (z^1(q), \ldots , z^m(q))$ is a solution to
(\ref{eq:three}).  Then $\phi$ is invariant under
$v=a^j(x)\frac{\partial}{\partial x^j}$ if and only if
\begin{equation}
\left| \begin{array}{cccc}
w^1(\alpha , z(q))  & \partial _2F^1 & \ldots
 & \partial _mF^1 \\
w^2(\alpha , z(q)) & \partial _2F^2 &  \ldots & \partial _mF^2 \\
\vdots & \vdots & \ddots & \vdots \\
w^m(\alpha , z(q)) & \partial _2F^m &  \ldots & \partial _mF^m
 \end{array}
\right| =0 \label{eq:red}
\end{equation}
for all $q\in U$, where $\alpha =
\alpha ^j\frac{\partial}{\partial q^j}$ is the complex vector
field associated to $v$ obtained by setting $\alpha
^j=a^{2j-1}+ia^{2j} \ (j=1, \ldots , m)$, and $w^i(q, z)\equiv
q^i- M^i_{\bar{j}}(z)q^{\bar{j}}$ (thus $w^i(q, z)$ is the part of $F^i$
which is homogeneous of degree $1$ in $q$).
\end{proposition}

\proof  The function $z^1$ is invariant under $v$ if and only if 
$$
\alpha ^I\frac{\partial z^1}{\partial q^I} = 0
$$
where we sum over $I=1, \ldots , m, \bar{1}, \ldots , \bar{m}$.
But 
$$
\alpha ^I\frac{\partial z^1}{\partial q^I}=-\alpha ^I(K^{-1})^1_b
\frac{\partial F^b}{\partial q^I}= -\alpha ^I\frac{1}
{\det K}\hat{K}^b_1\frac{\partial F^b}{\partial q^I} $$
where $\hat{K}^b_a$ is the $(b,a)$-th entry in the matrix of
cofactors of $K$.  But since the expression $ w^i(\alpha , z) =
\alpha ^i -M^i_{\bar{j}}(z)\alpha ^{\bar{j}}$ is homogeneous of
degree one in $\alpha$,
$$
\alpha ^I\frac{\partial F^b}{\partial q^I}=
\alpha ^I\frac{\partial w^b}{\partial q^I}=w^b(\alpha , z)
$$ and the result follows.
\eproof

Special cases of the above proposition are as follows:

(i) {\em Orthogonal projection\/}  $\pi^{2m}_1:{\Rmm}\rightarrow
{\RR}^{2m-1}$ given by the formula $\pi^{2m}_1(x^1,x^2, \ldots
, x^{2m}) = (x^2, \ldots , x^{2m})$.  Putting $a=(1, 0, \ldots ,
0)$, we see that $z^1$ reduces to ${\RR}^{2m-1}$ if and only if
\begin{equation}
\left| \begin{array}{cccc}
1  & \partial _2F^1 & \ldots
 & \partial _mF^1 \\
\mu _1 & \partial _2F^2 &  \ldots & \partial _mF^2 \\
\vdots & \vdots & \ddots & \vdots \\
\mu _{m-1} & \partial _2F^m &  \ldots & \partial _mF^m
 \end{array}
\right| =0\;. \label{eq:orth1}
\end{equation}

(ii) {\em Radial projection\/}  $\pi^{2m}_r:{\Rmm}\backslash\{
0\}\rightarrow S^{2m-1}$ given by the formula $\pi ^{2m}_r(x^1,
\ldots , x^n)=(x^1, \ldots , x^n)/|x|$. Put $a=(x^1, \ldots ,
x^{2m})$ so that $\alpha ^j=q^j$. In particular $\alpha
^I\displaystyle{\frac{\partial F^b}{\partial q^I}}=w^b(q, z)
=h^b(z)$, so that $z^1$ reduces to $S^{2m-1}$ if and only if
\begin{equation}
\left| \begin{array}{cccc}
h^1  & \partial _2F^1 & \ldots
 & \partial _mF^1 \\
h^2 & \partial _2F^2 &  \ldots & \partial _mF^2 \\
\vdots & \vdots & \ddots & \vdots \\
h^m & \partial _2F^m &  \ldots & \partial _mF^m
 \end{array}
\right| =0 \;. \label{eq:radial}
\end{equation}

(iii)  {\em Projection\/} $\pi ^{2m}_H:{\Rmm}\backslash\{
0\}\rightarrow {\CC}P^{m-1}$ given by $\pi^{2m}(q^1, \ldots ,
q^m) \mapsto [q^1, \ldots , q^m]$.  
The fibres of this map are spanned by the vector fields
${\displaystyle \alpha_1=
q^I\frac{\partial}{\partial q^I}}$ and
${\displaystyle \alpha_2=iq^i\frac{\partial}
{\partial q^i}-iq^{\bar{i}}\frac{\partial}{\partial
q^{\bar{i}}}}$ so that $z^1$ reduces to ${\CC}P^{m-1}$
if and only if (\ref{eq:radial}) above holds and
\begin{equation}
\left| \begin{array}{cccc}
q^1  & \partial _2F^1 & \ldots
 & \partial _mF^1 \\
q^2 & \partial _2F^2 &  \ldots & \partial _mF^2 \\
\vdots & \vdots & \ddots & \vdots \\
q^m & \partial _2F^m &  \ldots & \partial _mF^m
 \end{array}
\right| =0 \;. \label{eq:complexred}
\end{equation}

Finally, we wish to know when examples
are genuinely new examples and not simply obtained by the
composition with an orthogonal projection.  We therefore make the
definition:

\begin{definition} \cite{Ba-Wo-5}
Call a map $\phi:U\rightarrow \CC$ from an open subset of $\Rn$
{\em full\/} if we cannot write it as $\phi
=\psi\circ\pi _A$ for some orthogonal projection $\pi _A$ onto a
subspace $A$ of $\Rn$ and map $\psi :\pi _A(U)\rightarrow
{\CC}$.
\end{definition}

A test for fullness is given by Proposition 4.2 of
\cite{Ba-Wo-5}.   

\section{Heierarchies of Weierstrass representations}
\label{sec:weierstrass}
Let $J$ be a Hermitian structure defined on a domain $U$ of
${\Rmm}$ and suppose that $\phi :(U, J)\rightarrow {\CC}$ is
holomorphic.  As in Section \ref{sec:red}, in a neighbourhood of a
regular point we may assume that $\phi
=z^1$ where $z = (z^1, \ldots , z^m)$ is a solution $z(q)$ to Equation
(\ref{eq:three}):
$$
F(q,z) \equiv q-M(z)\bar{q}-h(z)=0
$$
for some holomorphic functions $\mu_i(z), (i=1,\ldots m(m-1)/2)$,
$h^i(z), (i=1,\ldots,m)$.
Throughout this section we assume that $z^1$ has superminimal
fibres (and so, in particular, is harmonic).  Then $M(z)$ depends
on $z^1$ only and the Jacobian matrix takes the form
$$
K=\left( \begin{array}{cccc}
\partial _1F^1  & \partial _2h^1 & \ldots
 & \partial _mh^1 \\
\partial _1F^2 & \partial _2h^2 &  \ldots & \partial _mh^2 \\
\vdots & \vdots & \ddots & \vdots \\
\partial _1F^m & \partial _2h^m &  \ldots & \partial _mh^m
 \end{array}
\right)\;.
$$

By assumption, the determinant of $K$ is non-zero on $U$, so that,
at least one of the cofactors $\hat{K}^1_1,\hat{K}^2_1, \ldots ,
\hat{K}^m_1$ is non-zero, say the cofactor $\hat{K}^a_1$ obtained
by omitting
row $a$ and column $1$ from $K$. Then, by the Implicit Function
Theorem, we can locally solve the equations $F^1=0, \ldots , F^{a-1}=0,
F^{a+1}=0, \ldots , F^m=0$ for $z^2, \ldots , z^m$ as holomorphic
functions of $z^1$, on substituting these functions into the remaining
equation $F^a=0$, this takes the form (\ref{eq:wei1}) in the following
proposition (cf. \cite[Proposition 3.12]{Ba-Wo-5}):

\begin{proposition}\label{prop:even}
Let $\mu _1(z^1), \ldots ,\mu _{m(m-1)/2}(z^1)$ be given
holomorphic functions of one variable and $\Psi _{2m}(w^1, \ldots
, w^{m+1})$ a given holomorphic function of $m+1$ variables. 
Consider the equation:
\begin{equation}
\widetilde{\Psi}_{2m}(q,z^1) \equiv \Psi _{2m}\left(
q^1-M^1_{\bar{j}}(z^1)q^{\bar{j}}, q^2-M^2_{\bar{j}}(z^1)
q^{\bar{j}}, \ldots , q^m-M^m_{\bar{j}}(z^1)q^{\bar{j}},
z^1\right) =0 \,. \label{eq:wei1}
\end{equation}
Suppose that, at a point $(q, z^1)$ satisfying $\widetilde{\Psi}_{2m}
(q, z^1)=0$,
$$
\frac{\partial\widetilde{\Psi}_{2m}}{\partial z^1} \neq 0 \;. 
$$
Then the local solution $z^1 = z^1(q)$ to  Equation (\ref{eq:wei1})
through that point is a
holomorphic harmonic morphism with superminimal fibres.  All
submersive holomorphic harmonic morphisms with superminimal
fibres are given this way locally.
\end{proposition}

We refer to Equation (\ref{eq:wei1}) as the {\em Weierstrass
representation for holomorphic harmonic morphisms with
superminimal fibres on\/} ${\Rmm}$.

\begin{example}\label{ex:4d}
If $m=2$, then Equation (\ref{eq:wei1}) has the form
$$
\Psi _4\left( q^1-\mu _1(z^1)q^{\bar{2}},
q^2+\mu _1(z^1)q^{\bar{1}}, z^1\right) = 0 \,,
$$
which is precisely the equation obtained in \cite{Wo} describing
{\em all\/} (submersive complex-valued) harmonic morphisms $z^1 = z^1(q)$
on domains of ${\RR}^4$.
\end{example}

\begin{example}\label{ex:tgeven}
Let $\alpha_1,\alpha_2,\ldots ,\alpha_m$ be $m$ holomorphic functions of
$z^1$ and consider the particular form of $\Psi _{2m}$ given by 
$$
\Psi _{2m}(w^1,\ldots , w^m, z^1)=\alpha_1w^1+\alpha _2w^2+\cdots +
\alpha_mw^m - 1\,.
$$
Then Equation (\ref{eq:wei1}) can be written
$$
(\alpha_1-M^k_{\bar{1}}\alpha_k)x^1+\ii (\alpha_1+M^k_{\bar{1}}\alpha_k)x^2
+ \cdots + (\alpha_m-M^k_{\bar{m}}\alpha_k)x^{2m-1}+\ii
(\alpha_m+M^k_{\bar{m}}\alpha_k)x^{2m}=1\,.
$$
Note that the sum of the squares of the coefficients of $x^1, \ldots ,x^{2m}$
equals $-4\sum M^k_{\bar{j}}\alpha_j\alpha_k$\,, which vanishes since
$M^k_{\bar{j}}$ is antisymmetric in $k$ and $j$.  We thus obtain the local
characterization (\ref{eq:tginhomog}) of those (submersive complex-valued)
harmonic morphisms with totally geodesic fibres,
on even-dimensional Euclidean spaces. 
\end{example}
 
\begin{example}\label{ex:6d}
(see also \cite{Ba-Wo-5}, Example 4.9)  Let $m=3$ and set $\Psi _6
(w^1, w^2, w^3, z^1)=w^1w^2-w^3, \mu _1=\mu _2=z^1, \mu _3=0$.
Then (\ref{eq:wei1}) becomes the quadratic equation
$$
(z^1)^2q^{\bar{1}}(q^{\bar{2}}+q^{\bar{3}})+(z^1)\left[
q^{\bar{1}}+q^2 (q^{\bar{2}}+q^{\bar{3}})-|q^1|^2\right]
+(q^3-q^1q^2) = 0 \,. $$
Any local solution is a full harmonic morphisms with superminimal
fibres. It is easy to see that $z^1$ is not holomorphic with
respect to any K\"ahler structure.

More generally, for arbitrary $m$, we can set $\Psi _{2m}(w^1, \ldots ,
w^m, z^1) =w^1w^2\ldots w^{m-1}-w^m$ to obtain generalizations to
arbitrary even dimensions.
\end{example}

Now suppose that $z^1$ factors through the projection $\pi
^{2m}_1: {\Rmm}\rightarrow {\RR}^{2m-1}$.  That is, from
Equation (\ref{eq:orth1}),
$$
\left| \begin{array}{cccc}
1  & \partial _2h^1 & \ldots
 & \partial _mh^1 \\
\mu _1 & \partial _2h^2 &  \ldots & \partial _mh^2 \\
\vdots & \vdots & \ddots & \vdots \\
\mu _{m-1} & \partial _2h^m &  \ldots & \partial _mh^m
 \end{array}
\right| =0\;.
$$  
After elementary row operations this becomes the $(m-1)\times
(m-1)$ determinant:
\begin{equation}
\left| \begin{array}{ccc}
\partial _2(h^2-\mu _1h^1) & \ldots & \partial _m(h^2-\mu _1h^1) \\
\partial _2(h^3-\mu _2h^1) & \ldots & \partial _m(h^3-\mu _2h^1) \\
\vdots & \ddots & \vdots \\
\partial _2(h^m-\mu _{m-1}h^1) & \ldots &
\partial _m(h^m-\mu _{m-1}h^1)
\end{array}\right| =0\;.\label{eq:det}
\end{equation}

We consider the special solution of (\ref{eq:det}) given by 
\begin{equation}
\partial _m(h^2-\mu _1h^1)=\partial _m(h^3-\mu _2h^1)=\ldots 
=\partial _m(h^m-\mu _{m-1}h^1)=0\;, \label{eq:sp1}
\end{equation}
that is 
\begin{equation}
\left\{ \begin{array}{ccc}
h^2-\mu _1h^1 & = & \alpha ^2(z^1, \ldots , z^{m-1}) \\
h^3-\mu _2h^1 & = & \alpha ^3(z^1, \ldots , z^{m-1}) \\
\vdots & \vdots & \vdots \\
h^m-\mu _{m-1}h^1 & = & \alpha ^m(z^1, \ldots , z^{m-1}) 
\end{array}\right. 
\label{eq:R}
\end{equation}
are all holomorphic functions independent of $z^m$.  Consider
the matrix
$$
A = \left( \partial _i\alpha ^j\right)
_{i=2, \ldots ,m-1}^{j=2, \ldots ,m} \,.
$$
Then amongst the minors $A^a=\det \left(\partial _i\alpha
^j\right) _{i=2,\ldots ,m-1}^{j=2, \ldots ,\hat{a}, \ldots , m}$
at least one is not identically zero otherwise $\det K$ would
vanish.  So amongst the $m-1$ equations (\ref{eq:R}), there are
$m-2$ for which we can eliminate $z^2, \ldots , z^{m-1}$. 
Substituting $h^i=q^i-M^i_{\bar{j}}(z^1)q^{\bar{j}}$
into the remaining equation gives the functional relation:
\begin{eqnarray*}
\lefteqn{\Phi _{2m-1}\left( q^2-q^{\bar{j}}M^2_{\bar{j}}(z^1)-
\mu _1(z^1)\big(q^1-q^{\bar{j}}
M^1_{\bar{j}}(z^1)\big), \ldots\right. }\\ 
& & \left. \ldots ,  q^m-q^{\bar{j}}M^m_{\bar{j}}(z^1)-
\mu _{m-1}(z^1)\big(q^1-q^{\bar{j}}M^1_{\bar{j}}(z^1)\big),
z^1\right) =0\;,
\end{eqnarray*}
where $\Phi _{2m-1}=\Phi _{2m-1}(u_1, \ldots , u_m)$ is a
holomorphic function of $m$ complex variables.

If now this equation is satisfied, the restriction of $z^1$ to any
level hypersurface $x^1=$ constant is a harmonic morphism by
Section \ref{sec:red}. For convenience we take $x^1=0$, and then,
setting $q^1=\ii x^2$, we obtain the following characterization:

\begin{proposition}\label{prop:odd}
Let $\mu _1(z^1), \ldots ,\mu _{m(m-1)/2}(z^1)$ be given
holomorphic functions of one variable and $\Phi _{2m-1}=\Phi
_{2m-1}(u^1, \ldots , u^m)$ a holomorphic function of $m$
variables.  Consider the equation
\pagebreak[3]
\begin{eqnarray}
\lefteqn{\widetilde{\Phi}_{2m-1}(q,z^1) \equiv
\Phi _{2m-1}\left( q^2  - 2\ii \mu _1(z^1)x^2
-\sum _{\bar{j}\geq 2}\left( M^2_{\bar{j}}(z^1)
-\mu _1(z^1) M^1_{\bar{j}}(z^1)\right) q^{\bar{j}}, \cdots\right. }
\nonumber \\
& & \left. \cdots, 
q^m-2\ii\mu _{m-1}(z^1)x^2
-\sum _{\bar{j}\geq 2}\left( M^m_{\bar{j}}(z^1)
-\mu _{m-1}(z^1) M^1_{\bar{j}}(z^1)\right) q^{\bar{j}}, z^1\right)
=0\;. \nonumber\\ & & \label{eq:wei2}
\end{eqnarray}

Suppose that, at a point $(q,z^1) = (x^2,q^2,\ldots q^m, z^1)$
satisfying $\widetilde{\Phi}_{2m-1}(q, z^1)=0$,
$$
\frac{\partial\widetilde{\Phi}_{2m-1}}{\partial z^1} \neq 0 \;. 
$$
Then the local solution $z^1 = z^1(q)$ to Equation (\ref{eq:wei2})
through that point is a
harmonic morphism $\phi :U\subset {\RR}^{2m-1}\rightarrow {\CC}$
whose lift\/ $\Phi =\phi \circ \pi ^{2m}_1$ to ${\Rmm}$
is a holomorphic harmonic morphism which has superminimal fibres and
satisfies the simplifying assumption (\ref{eq:sp1}).  All such
submersive harmonic morphisms on domains of ${\RR}^{2m-1}$ are
given this way locally.
\end{proposition}

\begin{example}\label{ex:R3}
If $m=2$, then the assumption (\ref{eq:sp1}) is automatic and 
(\ref{eq:wei2}) takes the form:
$$
\Phi_3(q^2-2\ii\mu _1(z^1) x^2-\big(\mu _1(z^1)\big)^2
q^{\bar{2}}, z^1)=0
$$
which is the local representation of {\em all\/} (submersive complex-valued)
harmonic morphisms $z^1 = z^1(x^2, q^2)$ on domains of
${\RR}^3$ \cite{Ba-Wo-1}.
\end{example}

\noindent {\bf Remark} \ By \cite[Lemma 4.3]{Ba-Wo-1} and \cite[Theorem
2.19]{Ba-Wo-2}, any harmonic morphism from
a domain of $\RR^3$
is locally a submersive complex-valued harmonic morphism followed by a
weakly conformal map.  In
higher dimensions, little is known about the behaviour of a harmonic
morphism near a critical
point: for some partial results in $4$ dimensions, see \cite{Wo}.
\medskip

\begin{example}\label{ex:tgodd}
Let $\alpha_1,\alpha_2,\ldots ,\alpha_{m-1}$ be $m-1$ holomorphic
functions of $z^1$ and consider the particular form of $\Phi _{2m-1}$
given by 
$$
\Phi _{2m-1}(w^1,\ldots , w^{m-1}, z^1)=\alpha_1w^1+\alpha _2w^2+\cdots
+\alpha_{m-1}w^{m-1} - 1\,.
$$ 
Then Equation (\ref{eq:wei2}) becomes
\begin{eqnarray*}
\left( -\sum ^{m-1}_{j=2}2\ii \alpha_j\mu_j\right) x^2  +  \left(
\alpha_1-\sum_k\alpha_k\left( M^k_{\bar{2}} -\mu
_kM^1_{\bar{2}}\right)\right)x^3 
  +  \ii\left(
\alpha_1+\sum_k\alpha_k\left( M^k_{\bar{2}} -\mu
_kM^1_{\bar{2}}\right)\right)x^4 + \cdots \\
\cdots  +  \left(
\alpha_{m-1}-\sum_k\alpha_k\left( M^k_{\bar{m}} -\mu
_kM^1_{\bar{m}}\right)\right)x^{2m-1} 
  +  \ii\left(
\alpha_1+\sum_k\alpha_k\left( M^k_{\bar{m}} -\mu
_kM^1_{\bar{m}}\right)\right)x^{2m}=1 \,.   
\end{eqnarray*}
As in Example \ref{ex:tgeven}, the sum of the squares of the
coefficients vanishes and we retrieve the local characterization
(\ref{eq:tginhomog}) of those (submersive complex-valued)
harmonic morphisms with totally geodesic fibres, now defined on
odd-dimensional Euclidean spaces.
\end{example}
 
\begin{example}\label{ex:R5}
Let $m=3$ and set $\Phi _5(w^1, w^2, z^1)=w^1w^2-1, \mu _1=\mu
_2=z^1, \mu _3 =0$.  Then (\ref{eq:wei2}) becomes the quartic
equation:
\begin{eqnarray*}
\lefteqn{(z^1)^4\big( q^{\bar{2}}+q^{\bar{3}}\big) ^2-4i(z^1)^3x^2
\big( q^{\bar{2}}+q^{\bar{3}}\big) } \\
& &  +(z^1)^2\big( q^2+q^3\big)\big( q^{\bar{2}}+q^{\bar{3}}
\big) -2\ii(z^1)x^2\big(q^2+q^3\big) +(q^2q^3-1)=0 \,.
\end{eqnarray*}
Any local solution $z^1=z^1(x^2, q^2, q^3)$ is a full
harmonic morphism defined on a domain of ${\RR}^5$.
\end{example}

Suppose that we now reduce once more by the vector field $\partial
/\partial x^2$. Then Condition (\ref{eq:red}) becomes:
\begin{equation}
\left| \begin{array}{cccc}
-1  & \partial _2h^1 & \ldots
 & \partial _mh^1 \\
\mu _1 & \partial _2h^2 &  \ldots & \partial _mh^2 \\
\vdots & \vdots & \ddots & \vdots \\
\mu _{m-1} & \partial _2h^m &  \ldots & \partial _m h^m
 \end{array}
\right| =0 \,, \label{eq:orth2}
\end{equation}
and (\ref{eq:orth1}) and (\ref{eq:orth2}) are satisfied
if and only if (\ref{eq:orth1}) holds and 
\begin{equation}
\left| \begin{array}{ccc}
\partial _2h^2 & \ldots & \partial _mh^2 \\
 \partial _2h^3 &  \ldots & \partial _mh^3 \\
\vdots & \ddots & \vdots \\
\partial _2h^m &  \ldots & \partial _mh^m
 \end{array}
\right| =0\,.
\label{eq:E}
\end{equation}

We consider the special case of (\ref{eq:E}) given by 
\begin{equation}
\partial _mh^2=\partial _mh^3=\ldots \partial _mh^m=0 \,.
\label{eq:sp2}
\end{equation}
Note, in particular, that (\ref{eq:sp2}) implies (\ref{eq:sp1}) and so 
such maps are a subset of those satisfying (\ref{eq:wei2}).  Then
$$
\left\{ 
\begin{array}{ccc}
h^2 & = & \beta ^2(z^1, \ldots , z^{m-1}) \\
h^3 & = & \beta ^3(z^1, \ldots , z^{m-1}) \\
\vdots & \ddots & \vdots \\
h^m & = & \beta ^m(z^1, \ldots , z^{m-1}) 
\end{array} \right.
$$
are all independent of $z^m$.  By eliminating $z^2, \ldots ,
z^{m-1}$ we now obtain a representation of $z^1$ in the form
$$
\Psi _{2m-2}\left( q^2-M^2_{\bar{j}}(z^1) q^{\bar{j}},q^3-
M^3_{\bar{j}}(z^1) q^{\bar{j}},
\ldots , q^m-M^m_{\bar{j}}(z^1) q^{\bar{j}}, z^1\right) =0
$$ 
where $\Psi _{2m-2}(w^2, \ldots , w^{m+2})$ is a holomorphic
function of $m$ complex variables.  But this is precisely
(\ref{eq:wei1}) with $m$ replaced by $m-1$.

\medskip

To sum up, we define data for the Weierstrass representation:

\medskip
{\bf Data}
For each $m=2, 3, \ldots $, let $\mu _1(z^1), \ldots ,\mu
_{m(m-1)/2}(z^1)$ be given holomorphic functions of one variable
and let $M=M(z^1)$ be the matrix given by (\ref{M}). Suppose that
either:

(i)  $\Psi _{2m}$ is a holomorphic function of $m+1$ complex
variables, or

(ii) $\Phi _{2m-1}$ is a holomorphic function of $m$ complex
variables.

\begin{theorem}\label{th:wei}
Suppose that we are given the data above.

(a) The equation
$$
\Psi _{2m}\left( q^1-M^1_{\bar{j}}(z^1) q^{\bar{j}},q^2-
M^2_{\bar{j}}(z^1)q^{\bar{j}},
\ldots , q^m-M^m_{\bar{j}}(z^1) q^{\bar{j}}, z^1\right) =0
$$
locally determines all submersive holomorphic harmonic morphisms
with superminimal fibres defined on domains of ${\Rmm}$ .

(b)  The equation
\begin{eqnarray*}
\lefteqn{\Phi _{2m-1}\left( q^2  - 2\ii\mu _1(z^1)x^2
-\sum _{\bar{j}\geq 2}\left( M^2_{\bar{j}}(z^1)
-\mu _1(z^1) M^1_{\bar{j}}(z^1)\right) q^{\bar{j}}, \ldots\right. }\\
& & \left. \ldots , 
q^m-2\ii\mu _{m-1}(z^1)x^2
-\sum _{\bar{j}\geq 2}\left( M^m_{\bar{j}}(z^1)
-\mu _{m-1}(z^1)M^1_{\bar{j}}(z^1)\right) q^{\bar{j}}, z^1\right) =0
\end{eqnarray*}
locally determines all submersive harmonic morphisms on domains of
${\RR}^{2m-1}$ 
which are the reduction of holomorphic harmonic morphisms on
domains of ${\Rmm}$ with superminimal fibres satisfying the
simplifying assumption (\ref{eq:sp1}).  Amongst the solutions are
the lifts of all submersive
holomorphic harmonic morphisms with superminimal
fibres defined on domains of ${\RR}^ {2m-2}$.

(c)  The equations $\Phi _3=0$ and $\Psi _4=0$ locally describe
{\em all\/} submersive harmonic morphisms on domains of ${\RR}^3$ and
${\RR}^4$, respectively.  
\end{theorem}

Schematically we represent the above heierarchy of representations
by the inclusions:
$$
\cdots \supset \{\Psi _{2m} =0 \} \supset \{\Phi _{2m-1} =0 \} \supset
\{\Psi_{2m-2} =0 \} \supset\cdots \supset \{\Psi _4 =0 \} \supset
\{\Phi _3 =0 \} \,.
$$

\section{Reduction to Euclidean spheres}\label{sec:spheres}
Let $U$ be an open subset of ${\Rmm}$ on which is defined a
Hermitian structure $J$ and let $\phi =z^1:U\rightarrow {\CC}$
be a harmonic morphism holomorphic with respect to $J$ and with
superminimal fibres.  In particular by Proposition
\ref{prop:even} we can suppose that $z^1$ is determined by
(\ref{eq:wei1}):
$$
\Psi (w^1, w^2, \ldots , w^m, z^1)=0
$$
where $w^i=q^i-M^i_{\bar{j}} q^{\bar{j}}$ and
$M^i_{\bar{j}} = M^i_{\bar{j}}(z^1)$.
Then $z^1$ is invariant under radial projection if and only if 
$$
q^I\frac{\partial z^1}{\partial q^I}=0
$$
where, as usual, we sum over $I=1,\ldots , m, \bar{1}, \ldots , \bar{m}$.
\begin{lemma}\label{lemma:homog}
With the data above, the equation
$$
q^I\frac{\partial z^1}{\partial q^I}=0 \,,
$$
is satisfied if and only if the condition 
\begin{equation}
w^i\frac{\partial\Psi}{\partial w^i}=0 \mbox{ whenever } \Psi =0
\label{eq:H}
\end{equation}
is satisfied.
\end{lemma}

\proof  Differentiating (\ref{eq:wei1}) with respect to $q^k$ yields
$$
\frac{\partial \Psi}{\partial w^i}\frac{\partial w^i}{\partial
q^k} + \frac{\partial \Psi}{\partial z^1}\frac{\partial
z^1}{\partial q^k} =0 \,,
$$
and, differentiating $w^i$,
$$
\frac{\partial w^i}{\partial q^k} = \delta ^i_k-
\frac{\partial M^i_{\bar{j}}}{\partial z^1}\frac{\partial
z^1}{\partial q^k} q^{\bar{j}} \,.
$$
Also differentiating (\ref{eq:wei1}) and $w^i$ with respect to $q^{\bar{k}}$
yields
$$
\frac{\partial \Psi}{\partial w^i}\frac{\partial w^i}{\partial
q^{\bar{k}}} + \frac{\partial \Psi}{\partial z^1}\frac{\partial
z^1}{\partial q^{\bar{k}}} =0
$$
and
$$
\frac{\partial w^i}{\partial q^{\bar{k}}} = -M^i_{\bar{k}}-
\frac{\partial M^i_{\bar{j}}}{\partial z^1}\frac{\partial z^1}
{\partial q^{\bar{k}}}q^{\bar{j}} \,.
$$
Combining the first and third equations gives                     
$$
\frac{\partial \Psi}{\partial w^i}q^K\frac{\partial w^i}{\partial
q^K} + \frac{\partial \Psi}{\partial z^1}q^K\frac{\partial
z^1}{\partial q^K} =0 \,,
$$
so that ${\displaystyle q^K\frac{\partial z^1}{\partial q^K}=0}$
if and only if ${\displaystyle q^K\frac{\partial w^i}{\partial
q^K}\frac{\partial \Psi} {\partial w^i}=0}$. But substituting the
expressions for ${\displaystyle \frac{\partial w^i}{\partial
q^K}}$ above gives  
$$
q^K\frac{\partial w^i}{\partial q^K}\frac{\partial \Psi}{\partial
w^i}= \left( q^i-M^i_{\bar{k}}q^{\bar{k}}\right)
\frac{\partial\Psi}{\partial w^i} -\left( \frac{\partial
M^i_{\bar{j}}}{\partial z^1}q^{\bar{j}} \frac
{\partial\Psi}{\partial w^i}\right) q^K\frac{\partial
z^1}{\partial q^K} \;,
$$
and the result follows.
\eproof

We therefore have

\begin{theorem}\label{th:H}
Let $\phi =z^1:U\rightarrow {\CC}$ be a
harmonic morphism from a domain of ${\Rmm}$
determined by (\ref{eq:wei1}).  Then $\phi$
reduces to $S^{2m-1}$ if and only if the condition (\ref{eq:H}):
$$
w^i\frac{\partial\Psi}{\partial w^i}=0 \mbox{ whenever } 
\Psi =0
$$
is satisfied.
\end{theorem}

Clearly Condition (\ref{eq:H}) is satisfied if $\Psi$ is
homogeneous in $(w^1, \ldots , w^m)$.  In fact a partial converse
holds:

\begin{proposition}\label{pr:H}
Let $\Psi$ be an irreducible polynomial in $m$ complex variables 
$w^1, \ldots , w^m$. Then Condition (\ref{eq:H})
implies that $\Psi$ is homogeneous in $w^1, \ldots , w^m$.
\end{proposition}

The proof follows by combining the following two lemmas:

\begin{lemma}\label{lemma:K}
Let $\Psi$ be an irreducible polynomial in $m$ complex variables 
$w^1, \ldots , w^m$\,. Then Condition (\ref{eq:H}) implies that
\begin{equation}
w^i\frac{\partial\Psi}{\partial w^i}=k\Psi
\label{eq:K}
\end{equation}
for some $k\in {\CC}$.
\end{lemma}

\proof  Since $\Psi$ is irreducible, by the Nullstellensatz (see,
e.g. \cite{Mu}), Condition (\ref{eq:H}) implies that
$$
w^i\frac{\partial\Psi}{\partial w^i}=\alpha\Psi
$$
for some polynomial $\alpha$.  But if $\Psi$ has degree $k_j$ in
$w^j, \ (j=1, \ldots , m)$, then so does ${\displaystyle
w^i\frac{\partial\Psi}{\partial w^i}}$ and therefore $\alpha$ is
constant.
\eproof

\begin{lemma}\label{lemma:anal}
Let $\Psi$ be an irreducible analytic function in $m$ complex
variables $w^1, \ldots , w^m$.  Then
$$
w^i\frac{\partial\Psi}{\partial w^i}=k\Psi \mbox{ for some } 
k\in {\CC}
$$
if and only if $\Psi$ is homogeneous.
\end{lemma}

\proof  For fixed $(w^1_0, \ldots , w^m_0)$, set $\Psi _t=
\Psi (tw^1_0, \ldots , tw^m_0) = \Psi(w^1, \ldots , w^m)$ where
$w^j=tw^j_0$. Then by the chain rule,
$$
\frac{d\Psi _t}{dt}=\frac{\partial\Psi}{\partial
w^i}\frac{dw^i}{dt}=w^i_0 \frac{\partial\Psi}{\partial w^i}
$$
so that
$$
t\frac{d\Psi _t}{dt}=w^i\frac{\partial\Psi}{\partial w^i}=k\Psi _t \,.
$$
Hence
$$
\frac{d\Psi _t}{\Psi _t}=\frac{kdt}{t} \,.
$$
Integrating from $t=1$ to $t=T$ gives:
$$
\ln\Psi _t-\ln\Psi _1=k\ln T
$$
so that $\Psi _t=\Psi _1T^k$, i.e.
$$
\Psi (tw^1_0, \ldots , tw^m_0)=t^k\Psi (tw^1, \ldots , tw^m) \,,
$$
showing that $\Psi$ is homogeneous of degree $k$.
\eproof

\noindent {\bf Remark} \ Lemma \ref{lemma:K} and the Proposition are false
for $\Psi$ an irreducible {\em analytic\/} function.  For example
$\Psi = e^{w_1}w_1$ is irreducible (since $e^{w_1}$ is a unit)
and
$$
w^1\frac{\partial\Psi}{\partial w^1}=(1+w^1)\Psi  \,,
$$
so that Condition (\ref{eq:H}) is satisfied but not (\ref{eq:K}),
and $\Psi$ is not homogeneous.

\bigskip

The above results give a method for constructing examples on
odd-dimensional spheres:

\begin{example} \label{ex:tgoddsphere}
Consider the homogeneous analogue of Example \ref{ex:tgeven}, thus, as in
that example, choose $m$ holomorphic functions
$\alpha_1, \ldots , \alpha_m$ and set 
$\Psi (w^1,\ldots , w^m, z^1)=\alpha_1w^1+\cdots +\alpha_mw^m$\,.
Then $\Psi =0$ determines locally all (submersive complex-valued)
harmonic morphisms on $S^{2m-1}$ with totally geodesic fibres.
\end{example}

If we now choose $\Psi$ to be an irreducible polynomial of degree $\geq 2$\,,
we obtain new examples of harmonic morphisms on odd-dimensional spheres: 

\begin{example}\label{ex:S5}
Let $\Psi(w^1, w^2, w^3, z^1)=(w^1)^2-(z^1)^2\left(
(w^2)^2+(w^3)^2\right) $. Then the corresponding harmonic morphism
$z^1=z^1(q)$ defined on a suitable domain of $S^5$ is given implicitly by
the equation
$$
\left( q^1-\mu _1q^{\bar{2}}-\mu _2q^{\bar{3}}\right) ^2=
(z^1)^2\left( \left( q^2+\mu _1q^{\bar{1}}-\mu
_3q^{\bar{3}}\right) ^2+ \left( q^3+\mu _2q^{\bar{1}}+\mu
_3q^{\bar{2}}\right) ^2\right)
$$   
where $\mu _1, \mu _2, \mu _3$ are arbitrary holomorphic functions
of $z^1$. The generic regular fibre of $z^1$ extends to a
compact minimal submanifold of $S^5$, which, after the
change of coordinates $X= q^1-\mu _1q^{\bar{2}}-\mu
_2q^{\bar{3}}, \ Y=z^1\big( q^2+\mu _1q^{\bar{1}}-\mu _3q^{\bar{3}}\big)
, \ Z=z^1\big( q^3+\mu _2q^{\bar{1}}+\mu _3q^{\bar{2}}\big)$\,, 
is the level set $F=0$
of the polynomial function $F:{\CC}^3\rightarrow {\CC}, \
F(X,Y,Z)=X^2-Y^2-Z^2$, holomorphic with respect to a K\"ahler
structure on ${\RR}^6$.  However, this K\"ahler structure
varies from fibre to fibre and the minimal submanifolds are not
the level sets of a function holomorphic with respect to a {\em
fixed\/} K\"ahler structure.  In particular we obtain a foliation
of an open set of $S^5$ by minimal codimension 2 submanifolds.
\end{example}

In a similar vein, we can construct examples on even-dimensional
spheres by choosing the function $\Phi$ in the Weierstrass
representation (\ref{eq:wei2}) to be homogeneous in the first
$m-1$ variables.  Firstly, choosing $\Phi$ to be linear we have

\begin{example}
Let $\alpha_1, \ldots , \alpha_{m-1}$ be $m-1$ holomorphic functions of $z^1$
and let $\Phi$ be given by $\Phi (w^1, \ldots , w^{m-1},
z^1)=\alpha_1w^1+\cdots +\alpha_{m-1}w^{m-1}$\,.  Then $\Phi =0$
determines locally all (submersive complex-valued)
harmonic morphisms on $S^{2m-2}$ with totally geodesic fibres.
\end{example}

Secondly, if we choose $\Phi$ to be an irreducible polynomial of degree
$\geq 2$\,, we obtain new examples of harmonic morphisms on
even-dimensional spheres, for example:

\begin{example}\label{ex:S6}
Let $m=4$ and set $\mu_1=\mu_2=\mu_3=z^1\,, \mu_4=\mu_5=\mu_6=(z^1)^2$\,.
Let $\Phi$ be given by 
$$
\Phi (w^1, w^2, w^3, z^1)=(w^1)^2+w^2w^3\,.
$$
Then the corresponding harmonic morphism $z^1=z^1(q)$ defined on a suitable
domain of $S^6$ is given implicitly by the quartic equation
\begin{eqnarray*}
(z^1)^4\left(
(q^{\bar{2}})^2+(2q^{\bar{2}}+
q^{\bar{3}})(2q^{\bar{2}}+2q^{\bar{3}}+q^{\bar{4}})\right)  +  (z^1)^3\left(
-4\ii x^2q^{\bar{2}}-6\ii x^2q^{\bar{3}} - 2\ii x^2q^{\bar{4}}\right)  \\
+ (z^1)^2\left( -8(x^2)^2+2|q^2|^2+2|q^3|^2+q^4(2q^{\bar{2}}+q^{\bar{3}}) + 
q^3(2q^{\bar{2}}+q^{\bar{4}})\right)  \\
+  2\ii z^1x^2(2q^2-q^3-q^4)+(q^2)^2+q^3q^4=0\,.
\end{eqnarray*}
\end{example}

Finally, to find examples with {\em non-superminimal\/} fibres
which reduce to $S^{2m-1}$ seems much harder; we have no general
theory but give one family of examples:

\begin{example}\label{ex:S7}
Let $m=4$ and let $\mu _1=\mu _1(z^1, z^2)$ be an arbitrary
holomorphic function of two variables.  Set $\mu _2=z^1, \mu _3=\mu
_4=0, \mu _5=(z^1)^p, \mu _6=(z^1)^q, h^3=(z^1)^rz^2, h^4=z^2$
with $\partial _4h^2=0$ and $\partial _4h^1\neq 0$, for integers
$p,q,r$ not yet specified.  Then $\triangle z^1=0$, $\det K$ is
not identically zero and Condition (\ref{eq:radial}) for
reduction to $S^7$ is satisfied.  Note that the fibres
of $z^1$ are superminimal if and only if $\partial _2\mu _1=0$.

The map $z^1$ is defined by the 3rd and 4th equations of the
system (\ref{eq:three}):
$$
\left\{
\begin{array}{ccc}
q^3 + \mu _2q^{\bar{1}} - \mu _4q^{\bar{2}} - \mu _6q^{\bar{4}}
-h^3 & = & 0 \\
q^4 + \mu _3q^{\bar{1}} + \mu _5q^{\bar{2}} + \mu _6q^{\bar{3}} -
h^4 & = & 0  \end{array} \right. \,;
$$
that is, 
\begin{equation}
(z^1)^{r+q}q^{\bar{3}}+(z^1)^{r+p}q^{\bar{2}}+(z^1)^{r}q^4+(z^1)^{q}q^{\bar{4}}
-z^1q^3=0 \,. \label{eq:S7}
\end{equation}
The test for fullness (that $z^1$ does not factor through any
orthogonal projection -- cf. \cite[Proposition 4.2]{Ba-Wo-5})
requires that the equation in the complex vector $\alpha\in {\CC}^4$:
$$
(z^1)^{r+q}\alpha ^{\bar{3}}+(z^1)^{r+p}\alpha
^{\bar{2}}+(z^1)^{r}\alpha ^4+ (z^1)^{q}\alpha
^{\bar{4}}-z^1\alpha ^3=0 $$
has only the trivial solution
$\alpha =(\alpha ^1,\alpha ^2,\alpha ^3,\alpha ^4)=0$.  In
general $p,q,r$ can be chosen so that this is the case, e.g.
$r=1, p\geq 1, q\geq p+2$ will suffice, giving a family of full
harmonic morphisms on domains of the sphere $S^7$.  The fibres
are totally geodesic and (\ref{eq:S7}) is of the form
(\ref{eq:tghomog}) of the Introduction.  Similar constructions
are possible for $m=5,6,\ldots $ \ . 
\end{example}

\end{document}